\documentstyle[11pt,amssymb,amsfonts]{article}

\begin{document}
\newcommand{\p}{\parallel }
\makeatletter \makeatother
\newtheorem{th}{Theorem}[section]
\newtheorem{lem}{Lemma}[section]
\newtheorem{de}{Definition}[section]
\newtheorem{rem}{Remark}[section]
\newtheorem{cor}{Corollary}[section]
\renewcommand{\theequation}{\thesection.\arabic {equation}}

\title{{\bf Curvature of multiply warped products with an affine connection}
}

\author{Yong Wang\\
{\scriptsize \it School of Mathematics and Statistics, Northeast Normal University,},\\{\scriptsize \it Changchun Jilin, 130024, China; E-mail:  wangy581@nenu.edu.cn}}

\date{}
\maketitle

\begin{abstract} In this paper, we
 study the Einstein multiply warped products with a semi-symmetric non-metric
 connection and the multiply warped products with a semi-symmetric non-metric
 connection with constant scalar curvature, we apply our results to
 generalized Robertson-Walker spacetimes with a semi-symmetric non-metric connection
 and generalized Kasner spacetimes with a semi-symmetric non-metric
 connection and find some new examples of Einstein affine
 manifolds and affine manifolds with constant scalar curvature.
  We also consider the multiply warped products with an affine connection with a zero torsion. \\

\end{abstract}

\section{Introduction}
    \quad  The (singly) warped product $B\times_bF$  of two pseudo-Riemannian manifolds $(B,g_B)$ and
    $(F,g_F)$ with a smooth function $b:B\rightarrow (0,\infty)$ is the product
    manifold $B\times F$ with the metric tensor $g=g_B\oplus
    b^2g_F.$ Here, $(B,g_B)$ is called the base manifold and
    $(F,g_F)$ is called as the fiber manifold and $b$ is called as
    the warping function. Generalized Robertson-Walker space-times
    and standard static space-times are two well-known warped
    product spaces.The concept of warped products was first introduced by
    Bishop and ONeil (see [BO]) to construct examples of Riemannian
    manifolds with negative curvature. In Riemannian geometry,
    warped product manifolds and their generic forms have been used
    to construct new examples with interesting curvature properties
    since then. In [DD], F. Dobarro and E. Dozo had studied from the viewpoint of partial differential equations and variational methods,
    the problem of showing when a Riemannian metric of constant scalar curvature can be produced on a product manifolds by a warped product construction.
    In [EJK], Ehrlich, Jung and Kim got explicit solutions to warping function to have a constant scalar curvature for generalized Robertson-Walker space-times.
    In [ARS], explicit solutions were also obtained for the warping
    function to make the space-time as Einstein when the fiber is
    also Einstein.\\
      \indent One can generalize singly warped products to multiply warped
      products. Briefly, a multiply warped product $(M,g)$ is a
      product manifold of form
      $M=B\times_{b_1}F_1\times_{b_2}F_2\cdots \times_{b_m}F_m$ with
      the metric $g=g_B\oplus b_1^2g_{F_1}\oplus
      b_2^2g_{F_2}\cdots \oplus b_m^2g_{F_m}$, where for each $i\in
      \{ 1,\cdots,m\},~b_i:~B\rightarrow (0,\infty)$ is smooth and
      $(F_i,g_{F_i})$ is a pseudo-Riemannian manifold. In
      particular, when $B=(c,d)$ with the negative definite metric
      $g_B=-dt^2$ and $(F_i,g_{F_i})$ is a Riemannian manifold, we
      call $M$ as the multiply generalized Robertson-Walker
      space-time.  In
      [DU], Dobarro and \"{U}nal studied Ricci-flat and
      Einstein-Lorentzian multiply warped products and considered
      the case of having constant scalar curvature for multiply warped
      products and applied their results to generalized Kasner
      space-times.\\
      \indent Singly warped products have a natural generalization. A twisted product
      $(M,g)$ is a
      product manifold of form
      $M=B\times_{b}F$, with a smooth function $b:B\times F\rightarrow (0,\infty)$,
      and the metric tensor $g=g_B\oplus b^2g_F.$ In [FGKU], they
      showed that mixed Ricci-flat twisted products could be expressed as warped
      products. As a consequence, any Einstein twisted products are
      warped products. In [Wa], we define the multiply twisted
      products as generalizations of multiply warped
      products and twisted products. A multiply twisted product $(M,g)$ is a
      product manifold of form
      $M=B\times_{b_1}F_1\times_{b_2}F_2\cdots \times_{b_m}F_m$ with
      the metric $g=g_B\oplus b_1^2g_{F_1}\oplus
      b_2^2g_{F_2}\cdots \oplus b_m^2g_{F_m}$, where for each $i\in
      \{ 1,\cdots,m\},~b_i:~B\times F_i\rightarrow (0,\infty)$ is
      smooth.\\
              \indent The definition of a semi-symmetric metric connection was given by H. Hayden in [Ha]. In 1970, K. Yano [Ya]
        considered a semi-symmetric metric connection and studied some of its properties. He proved that a Riemannian manifold admitting
         the semi-symmetric metric connection has vanishing curvature tensor if and only if it is conformally flat. Motivated by the Yano' result,
         in [SO1], Sular and \"{O}zgur studied warped product manifolds with a
         semi-symmetric metric connection, they computed curvature of semi-symmetric metric connection
          and considered Einstein warped product manifolds with a semi-symmetric metric connection. In [Wa],
           we considered multiply twisted products with a semi-symmetric metric connection and
            computed the curvature of a semi-symmetric metric connection.
            We
         showed that mixed Ricci-flat multiply twisted products with a semi-symmetric metric
 connection can be expressed as multiply warped products which
 generalizes the result in [FGKU]. We also
 studied the Einstein multiply warped products with a semi-symmetric metric
 connection and multiply warped products with a semi-symmetric metric
 connection with constant scalar curvature, we applied our results to
 generalized Robertson-Walker spacetimes with a semi-symmetric metric connection
 and generalized Kasner spacetimes with a semi-symmetric metric
 connection and we found some new examples of Einstein affine
 manifolds and affine manifolds with constant scalar curvature. We
 also classified generalized Einstein Robertson-Walker spacetimes with a semi-symmetric metric
 connection and generalized Einstein Kasner spacetimes with a semi-symmetric metric
 connection. In [AC1,2], Agashe and Chafle introduced the notation
 of a semi-symmetric metric connection and studied some of its
 properties and submanifolds of a Riemannian manifold with semi-symmetric non-metric
 connections. In [SO2], they studied warped product manifolds with a
 semi-symmetric non-metric connection. The purpose of this paper is
 to study the Einstein multiply warped products with a semi-symmetric non-metric
 connection and multiply warped products with a semi-symmetric non-metric
 connection with constant scalar curvature. \\
      \indent This paper is arranged as follows: In Section 2,  we compute curvature
      of multiply twisted products with a semi-symmetric non-metric
 connection. In Section 3, we study the special multiply warped products with a
semi-symmetric non-metric connection.
  In Section 4, we study the generalized Robertson-Walker spacetimes with a semi-symmetric non-metric connection. In Section 5,
  we consider the generalized Kasner spacetimes with a semi-symmetric non-metric
 connection. In Section 6, we compute curvature
      of multiply twisted products with an affine
 connection with a zero torsion.
 \\

\section{Preliminaries}

 \noindent {\bf Definition 2.1}  A {\it multiply twisted product}
$(M,g)$ is a product manifold of form
      $M=B\times_{b_1}F_1\times_{b_2}F_2\cdots \times_{b_m}F_m$ with
      the metric $g=g_B\oplus b_1^2g_{F_1}\oplus
      b_2^2g_{F_2}\cdots \oplus b_m^2g_{F_m}$, where for each $i\in
      \{ 1,\cdots,m\},~b_i:~B\times F_i\rightarrow (0,\infty)$ is
      smooth.\\
\indent Here, $(B,g_B)$ is called the base manifold and
    $(F_i,g_{F_i})$ is called as the fiber manifold and $b_i$ is called as
    the twisted function. Obviously, twisted products and multiply warped
    products are the special cases of multiply twisted products.\\
\indent Let $M$ be a Riemannian manifold with Riemannian metric $g$.
A linear connection $\overline{\nabla}$ on a Riemannian manifold $M$
is called a {\it semi-symmetric connection} if the torsion tensor
$T$ of the connection $\overline{\nabla}$
$$T(X,Y)=\overline{\nabla}_XY-\overline{\nabla}_YX-[X,Y]\eqno(2.1)$$
satisfies
$$T(X,Y)=\pi(Y)X-\pi(X)Y,\eqno(2.2)$$
where $\pi$ is a 1-form associated with the vector field $P$ on $M$
defined by $\pi(X)=g(X.P).$ $\overline{\nabla}$ is called a {\it
semi-symmetric metric connection} if it satisfies
$\overline{\nabla}g=0.$  $\overline{\nabla}$ is called a {\it
semi-symmetric non-metric connection} if it satisfies
$\overline{\nabla}g\neq 0.$

 If $\nabla$ is the Levi-Civita connection of
$M$, the semi-symmetric non-metric connection $\overline{\nabla}$ is
given by
$$\overline{\nabla}_XY=\nabla_XY+\pi(Y)X,\eqno(2.3)$$
(see [AC1]). Let $R$ and $\overline{R}$ be the curvature tensors of
$\nabla$ and $\overline{\nabla}$ respectively. Then $R$ and
$\overline{R}$ are related by
$$\overline{R}(X,Y)Z=R(X,Y)Z+g(Z,\nabla_XP)Y-g(Z,\nabla_YP)X+
\pi(Z)[\pi(Y)X-\pi(X)Y],\eqno(2.4)$$ for any vector fields $X,Y,Z$
on $M$ [AC1]. By (2.3) and Proposition
2.2 in [Wa], we have\\

 \noindent {\bf Proposition 2.2} {\it Let $M=B\times_{b_1}F_1\times_{b_2}F_2\cdots
    \times_{b_m}F_m$ be a multiply twisted product and let $X,Y\in \Gamma(TB)$
      and $U\in \Gamma(TF_i)$, $W\in \Gamma(TF_j)$ and $P\in \Gamma(TB)$ . Then}\\
      \noindent$(1) ~~\overline{\nabla}_XY=\overline{\nabla}^B_XY.$\\
      \noindent$(2)~~\overline{\nabla}_XU=\frac{X(b_i)}{b_i}U.$\\
\noindent$(3)~~\overline{\nabla}_UX=[\frac{X(b_i)}{b_i}+\pi(X)]U.$\\
            \noindent$(4)~~\overline{\nabla}_UW=0~~if~ i\neq j.$\\
      \noindent$(5)~~\overline{\nabla}_UW=U({\rm ln}b_i)W+W({\rm ln}b_i)U-\frac{g_{F_i}(U,W)}{b_i}{\rm
      grad}_{F_i}b_i-b_ig_{F_i}(U,W){\rm
      grad}_{B}b_i+\nabla^{F_i}_UW~~if~ i= j.$\\

\noindent {\bf Proposition 2.3} {\it Let
$M=B\times_{b_1}F_1\times_{b_2}F_2\cdots
    \times_{b_m}F_m$ be a multiply twisted product and let $X,Y\in \Gamma(TB)$
      and $U\in \Gamma(TF_i)$, $W\in \Gamma(TF_j)$ and $P\in \Gamma(TF_k)$ . Then}\\
      \noindent$(1) ~~\overline{\nabla}_XY={\nabla}^B_XY.$\\
      \noindent$(2)~~\overline{\nabla}_XU=\frac{X(b_i)}{b_i}U+g(P,U)X.$\\
\noindent$(3)~~\overline{\nabla}_UX=\frac{X(b_i)}{b_i}U.$\\
            \noindent$(4)~~\overline{\nabla}_UW=g(W,P)U~~if~ i\neq j.$\\
      \noindent$(5)~~\overline{\nabla}_UW=U({\rm ln}b_i)W+W({\rm ln}b_i)U-\frac{g_{F_i}(U,W)}{b_i}{\rm
      grad}_{F_i}b_i-b_ig_{F_i}(U,W){\rm
      grad}_{B}b_i+\nabla^{F_i}_UW+\pi(W)U~~if~ i= j.$\\

\indent By (2.4) and Proposition
2.4 in [Wa], we have\\

\noindent {\bf Proposition 2.4} {\it Let
$M=B\times_{b_1}F_1\times_{b_2}F_2\cdots
    \times_{b_m}F_m$ be a multiply twisted product and let $X,Y,Z\in \Gamma(TB)$
      and $V\in \Gamma(TF_i)$, $W\in \Gamma(TF_j)$, $U\in \Gamma(TF_k)$ and $P\in \Gamma(TB)$. Then}\\
\noindent $(1)\overline{R}(X,Y)Z=\overline{R}^B(X,Y)Z.$\\
\noindent $(2)\overline{R}(V,X)Y=-\left[\frac{H^{b_i}_B(X,Y)}{b_i}
+g(Y,\nabla_XP)-\pi(X)\pi(Y)\right]
V.$\\
\noindent $(3)\overline{R}(X,V)W=\overline{R}(V,W)X=\overline{R}(V,X)W=0 ~if ~i\neq j.$\\
\noindent $(4)\overline{R}(X,Y)V=0.$\\
\noindent $(5)\overline{R}(V,W)X=VX({\rm ln}b_i)W-WX({\rm ln}b_i)V ~ if ~i=j.$\\
\noindent $(6)\overline{R}(V,W)U=0~ if ~i=j\neq k~ or~ i\neq j \neq k.$\\
\noindent $(7)\overline{R}(U,V)W=-g(V,W)\frac{g_B({\rm
grad}_Bb_i,{\rm grad}_Bb_k)}{b_ib_k}U-g(V,W)\frac{P(b_i)}{b_i}U$
, ~if $~i=j\neq k.$\\
\noindent $(8)\overline{R}(X,V)W=[WX({\rm ln}b_i)]V-g(W,V)
\cdot\left[\frac{\nabla_X^B({\rm grad}_Bb_i)}{b_i}+\frac{{\rm
grad}_{F_i}(X{\rm ln}b_i)}{b_i^2}+\frac{P(b_i)}{b_i}X\right]
 ~if
~i=j.$\\
\noindent $(9)\overline{R}(U,V)W=g(U,W){\rm grad}_B(V({\rm
ln}b_i))-g(V,W){\rm grad}_B(U({\rm
ln}b_i))+R^{F_i}(U,V)W-\left(\frac{|{\rm
grad}_Bb_i|^2_B}{b_i^2}+\frac{P(b_i)}{b_i}\right)\left[g(V,W)U-g(U,W)V\right] ~if ~i=j=k.~$\\

\noindent {\bf Proposition 2.5} {\it Let
$M=B\times_{b_1}F_1\times_{b_2}F_2\cdots
    \times_{b_m}F_m$ be a multiply twisted product and let $X,Y,Z\in \Gamma(TB)$
      and $V\in \Gamma(TF_i)$, $W\in \Gamma(TF_j)$, $U\in \Gamma(TF_k)$ and $P\in \Gamma(TF_l)$. Then}\\
\noindent $(1)\overline{R}(X,Y)Z={R}^B(X,Y)Z
.$\\
\noindent $(2)\overline{R}(V,X)Y=-\frac{H^{b_i}_B(X,Y)}{b_i}V
~if ~i\neq l$\\
\noindent $(3)\overline{R}(V,X)Y=-\frac{H^{b_i}_B(X,Y)}{b_i}V
-\pi(V)\frac{Y(b_i)}{b_i}X
~if ~i= l$\\
\noindent $(4)\overline{R}(X,V)W=\frac{X(b_l)}{b_l}\pi(W)V~if ~i\neq j.$\\
\noindent
$(5)\overline{R}(V,W)X=-\delta^l_i\frac{\pi(V)}{b_i}X(b_i)W+\delta^l_j\frac{\pi(W)}{b_j}X(b_j)V~if ~i\neq j.$\\
\noindent $(6)\overline{R}(X,Y)V=\pi(V)\left[\frac{X(b_l)}{b_l}Y-\frac{Y(b_l)}{b_l}X\right].$\\
\noindent $(7)\overline{R}(V,W)X=VX({\rm ln}b_i)W-WX({\rm ln}b_i)V-\delta_i^l\frac{X(b_i)}{b_i}[\pi(V)W-\pi(W)V] ~ if ~i=j.$\\
\noindent $(8)\overline{R}(V,W)U=0~ if ~i=j\neq k~ or~ i\neq j \neq k.$\\
\noindent $(9)\overline{R}(U,V)W=-g(V,W)\frac{g_B({\rm
grad}_Bb_i,{\rm grad}_Bb_k)}{b_ib_k}U -g(W,\nabla_VP)U
+\pi(W)[\pi(V)U-\pi(U)V]
, ~if ~i=j\neq k.$\\
\noindent $(10)\overline{R}(X,V)W=[WX({\rm
ln}b_i)]V-g(W,V)\frac{\nabla_X^B({\rm grad}_Bb_i)}{b_i}-{\rm
grad}_{F_i}(X{\rm ln}b_i)g_{F_i}(W,V)
+\frac{X(b_l)}{b_l}\pi(W)V-g(W,\nabla_VP)X+\pi(V)\pi(W)X
 ~if
~i=j.$\\
\noindent $(11)\overline{R}(U,V)W=g(U,W){\rm grad}_B(V({\rm
ln}b_i))-g(V,W){\rm grad}_B(U({\rm
ln}b_i))+R^{F_i}(U,V)W-\frac{|{\rm
grad}_Bb_i|^2_B}{b_i^2}(g(V,W)U-g(U,W)V)~if ~i=j=k\neq l.~$\\
\noindent $(12)\overline{R}(U,V)W=g(U,W){\rm grad}_B(V({\rm
ln}b_i))-g(V,W){\rm grad}_B(U({\rm
ln}b_i))+R^{F_i}(U,V)W-\frac{|{\rm
grad}_Bb_i|^2_B}{b_i^2}(g(V,W)U-g(U,W)V)
+g(W,\nabla_UP)V-g(W,\nabla_VP)U +\pi(W)[\pi(V)U-\pi(U)V]
 ~if ~i=j=k=l.~$\\

By proposition 2.4 and 2.5, we have\\

\noindent {\bf Proposition 2.6} {\it Let
$M=B\times_{b_1}F_1\times_{b_2}F_2\cdots
    \times_{b_m}F_m$ be a multiply twisted product and let $X,Y,Z\in \Gamma(TB)$
      and $V\in \Gamma(TF_i)$, $W\in \Gamma(TF_j)$ and $P\in \Gamma(TB)$. Then}\\
\noindent $(1) \overline{{\rm Ric}} (X,Y)=\overline{{\rm
Ric}}^B(X,Y)+\sum_{i=1}^ml_i\left[ \frac{H_B^{b_i}(X,Y)}{b_i}
+g(Y,\nabla_XP)-\pi(X)\pi(Y)\right]
.$\\
\noindent $(2) \overline{{\rm Ric}} (X,V)=\overline{{\rm Ric}}
(V,X)=(l_i-1)[VX(lnb_i)].$\\
\noindent $(3) \overline{{\rm Ric}} (V,W)=0~if ~i\neq j.$\\
\noindent $(4) \overline{{\rm Ric}} (V,W)={\rm Ric}^{F_i}
(V,W)+\left[\frac{\triangle_Bb_i}{b_i}+(l_i-1)\frac{|{\rm
grad}_Bb_i|^2_B}{b_i^2}+\sum_{j\neq i}l_j\frac{g_B({\rm
grad}_Bb_i,{\rm grad}_Bb_j)}{b_ib_j}\right.$\\
$\left.+\sum_{j=
1}^ml_j\frac{Pb_j}{b_j}+({n}-1)\frac{Pb_i}{b_i}\right]
g(V,W)~if ~i= j,$\\
{\it where ${\rm dim}B=n,~{\rm
dim}M=\overline{n}.$}\\

\noindent {\bf Corollary 2.7} {\it Let
$M=B\times_{b_1}F_1\times_{b_2}F_2\cdots
    \times_{b_m}F_m$ be a multiply twisted product and ${\rm
    dim}F_i>1$ and $P\in \Gamma(TB)$, then $(M,\overline{\nabla})$ is mixed Ricci-flat if and only if $M$ can
    be expressed as a multiply warped product. In particular, if $(M,\overline{\nabla})$
    is Einstein, then $M$ can be expressed as a multiply warped product.}\\

\noindent {\bf Proposition 2.8} {\it Let
$M=B\times_{b_1}F_1\times_{b_2}F_2\cdots
    \times_{b_m}F_m$ be a multiply twisted product and let $X,Y,Z\in \Gamma(TB)$
      and $V\in \Gamma(TF_i)$, $W\in \Gamma(TF_j)$ and $P\in \Gamma(TF_r)$. Then}\\
\noindent $(1) \overline{{\rm Ric}} (X,Y)={{\rm
Ric}}^B(X,Y)+\sum_{i=1}^ml_i\frac{H_B^{b_i}(X,Y)}{b_i}
.$\\
\noindent $(2) \overline{{\rm Ric}} (X,V)=(l_i-1)[VX({\rm ln}b_i)]+(\overline{n}-1)\frac{X(b_r)}{b_r}\pi(V).$\\
\noindent $(3) \overline{{\rm Ric}} (V,X)=(l_i-1)[VX({\rm ln}b_i)]+(1-\overline{n})\frac{X(b_r)}{b_r}\pi(V).$\\
\noindent $(3) \overline{{\rm Ric}} (V,W)=0~if ~i\neq j.$\\
\noindent $(4) \overline{{\rm Ric}} (V,W)={\rm Ric}^{F_i}
(V,W)+g(V,W)\left[\frac{\triangle_Bb_i}{b_i}+(l_i-1)\frac{|{\rm
grad}_Bb_i|^2_B}{b_i^2}+\sum_{j\neq i}l_j\frac{g_B({\rm
grad}_Bb_i,{\rm
grad}_Bb_j)}{b_ib_j}\right]+(\overline{n}-1)g(W,\nabla_VP)+(1-\overline{n})\pi(V)\pi(W)~if
~i= j
,$\\

\noindent {\bf Corollary 2.9} {\it Let
$M=B\times_{b_1}F_1\times_{b_2}F_2\cdots
    \times_{b_m}F_m$ be a multiply twisted product and ${\rm
    dim}F_i>1$ and $P\in \Gamma(TF_r)$, then $(M,\overline{\nabla})$ is mixed Ricci-flat if and only if $M$ can
    be expressed as a multiply warped product and $b_r$ is only dependent on $F_r$. In particular, if $(M,\overline{\nabla})$
    is Einstein, then $M$ can be expressed as a multiply warped product.}\\

\noindent {\bf Proposition 2.10} {\it Let
$M=B\times_{b_1}F_1\times_{b_2}F_2\cdots
    \times_{b_m}F_m$ be a multiply twisted product and $P\in \Gamma(TB)$, then the scalar
    curvature $\overline{S}$ has the following expression:}\\
    $$
    \overline{S}=\overline{S}^B+2\sum_{i=1}^m\frac{l_i}{b_i}\triangle_Bb_i+\sum_{i=1}^m\frac{S^{F_i}}{b_i^2}+\sum_{i=1}^ml_i(l_i-1)
\frac{|{\rm grad}_Bb_i|^2_B}{b_i^2}$$
 $$+\sum_{i=1}^m\sum_{j\neq
i}l_il_j\frac{g_B({\rm grad}_Bb_i,{\rm grad}_Bb_j)}{b_ib_j}+
(n-1)\sum_{i=1}^ml_i\frac{P(b_i)}{b_i}$$
$$+\sum_{i=1}^m\sum_{j=
1}^ml_il_j\frac{P(b_j)}{b_j}+\sum_{i=1}^ml_i\left[{\rm
div}_BP-\pi(P)\right].\eqno(2.5)$$\\

\noindent {\bf Proposition 2.11} {\it Let
$M=B\times_{b_1}F_1\times_{b_2}F_2\cdots
    \times_{b_m}F_m$ be a multiply twisted product and $P\in \Gamma(TF_r)$, then the scalar
    curvature $\overline{S}$ has the following expression:}\\
    $$
    \overline{S}={S}^B+2\sum_{i=1}^m\frac{l_i}{b_i}\triangle_Bb_i+\sum_{i=1}^m\frac{S^{F_i}}{b_i^2}+\sum_{i=1}^ml_i(l_i-1)
\frac{|{\rm grad}_Bb_i|^2_B}{b_i^2}$$ $$+\sum_{i=1}^m\sum_{j\neq
i}l_il_j\frac{g_B({\rm grad}_Bb_i,{\rm grad}_Bb_j)}{b_ib_j}+
(1-\overline{n})\pi(P)+(\overline{n}-1)\sum_{j_r=1}^{l_r}\varepsilon_{j_r}g(\nabla_{E_{j_r}^r}P,E_{j_r}^r).\eqno(2.6)$$\\

\section{ Special multiply warped products with a
semi-symmetric non-metric connection}

 \quad Let $M=I\times_{b_1}F_1\times_{b_2}F_2\cdots
    \times_{b_m}F_m$ be a multiply warped product with the metric tensor $-dt^2\oplus  b_1^2g_{F_1}\oplus \cdots\oplus  b_m^2g_{F_m}$ and $I$ is an open
    interval in $\mathbb{R}$ and $b_i\in C^{\infty}(I)$.\\

\noindent {\bf Theorem 3.1} {\it  Let
$M=I\times_{b_1}F_1\times_{b_2}F_2\cdots
    \times_{b_m}F_m$ be a multiply warped product with the metric tensor $-dt^2\oplus b_1^2g_{F_1}\oplus \cdots\oplus
    b_m^2g_{F_m}$ and $P=\frac{\partial}{\partial t}$. Then
    $(M,\overline{\nabla})$ is Einstein with the Einstein constant
    $\lambda$ if and only if the following conditions are satisfied
    for any $i\in\{1,\cdots,m\}$}\\
\noindent (1){\it  $(F_i,\nabla^{F_i})$ is Einstein with the
Einstein
constant $\lambda_i$, $i\in\{1,\cdots,m\}$.}\\
\noindent (2)
$\sum_{i=1}^ml_i\left(1-\frac{b_i''}{b_i}\right)=\lambda.$\\
\noindent (3)$\lambda_i-b_ib_i''-(l_i-1)b_i'^2-b_ib_i'\sum_{j\neq
i}l_j\frac{b_j'}{b_j}+b_i^2\sum_{j= 1}^ml_j\frac{b_j'}{b_j}=\lambda
b_i^2.$\\

\noindent{\bf Proof.} By Proposition 2.6, we have
$$\overline{{\rm Ric}}\left(\frac{\partial}{\partial t},\frac{\partial}{\partial
t}\right)=-\sum_{i=1}^ml_i\left(1-\frac{b_i''}{b_i}\right);\eqno(3.1)$$
$$\overline{{\rm Ric}}\left(\frac{\partial}{\partial t},V\right)=\overline{{\rm Ric}}\left(V,\frac{\partial}{\partial
t}\right)=0;\eqno(3.2)$$
$$\overline{{\rm Ric}}\left(V,W\right)={\rm
Ric}^{F_i}(V,W)+g_{F_i}(V,W)\left[-b_ib_i''-(l_i-1)b_i'^2-b_ib_i'\sum_{j\neq
i}l_j\frac{b_j'}{b_j}+b_i^2\sum_{j=1}^ml_j\frac{b_j'}{b_j}\right].\eqno(3.3)$$
By (3.1)-(3.3) and the Einstein condition, we get the above
theorem.~~~~ $\Box$\\

 \noindent {\bf Definition 3.2} $(M,\overline{\nabla})$ is called
 {\it pseudo-Einstein} with the Einstein constant $\lambda$ if $\frac{1}{2}[\overline{{\rm
 Ric}}(X,Y)+\overline{{\rm
 Ric}}(Y,X)]=\lambda g(X,Y).$\\

\noindent {\bf Theorem 3.3} {\it  Let
$M=I\times_{b_1}F_1\times_{b_2}F_2\cdots
    \times_{b_m}F_m$ be a multiply warped product with the metric tensor $-dt^2\oplus b_1^2g_{F_1}\oplus \cdots\oplus
    b_m^2g_{F_m}$ and $P\in \Gamma(TF_r)$ and $\overline{n}>2$. Then
    $(M,\overline{\nabla})$ is pseudo-Einstein with the Einstein constant
    $\lambda$ if and only if the following conditions are satisfied
    for any $i\in\{1,\cdots,m\}$}\\
\noindent (1){\it  $(F_i,\nabla^{F_i})~(i\neq r)$ is Einstein with
the Einstein
constant $\lambda_i$, $i\in\{1,\cdots,m\}$.}\\
\noindent (2){\it
$-\sum_{i=1}^ml_i\frac{b_i''}{b_i}=\lambda.$}\\
\noindent (3){\it $ {\rm
Ric}^{F_i}(V,W)-g_{F_i}(V,W)\left[b_ib_i''+(l_i-1)b_i'^2+b_ib_i'\sum_{j\neq
i}l_j\frac{b_j'}{b_j}+\lambda b_i^2\right]$\\
$=(\overline{n}-1)\left[\pi(V)\pi(W)-\frac{g(W,\nabla_VP)+g(V,\nabla_WP)}{2}\right],~for~
V,W\in \Gamma(TF_r), ~r=i.$}\\
 \noindent
(4){\it $\lambda_i-b_ib_i''-(l_i-1)b_i'^2-b_ib_i'\sum_{j\neq
i}l_j\frac{b_j'}{b_j}-\lambda
b_i^2=0,$ for $i\neq r.$}\\

\noindent{\bf Proof.} By Proposition 2.8, then ,$$\overline{{\rm
Ric}}\left(\frac{\partial}{\partial t},\frac{\partial}{\partial
t}\right)=\sum_{i=1}^ml_i\frac{b_i''}{b_i};\eqno(3.4)$$ So we have
$\sum_{i=1}^ml_i\frac{b_i''}{b_i}=-\lambda.$

$$  \overline{{\rm Ric}} (V,W)={\rm Ric}^{F_i}
(V,W)+b_i^2g_{F_i}(V,W)\left[-\frac{b_i''}{b_i}+(l_i-1)\frac{-b_i'^2}{b_i^2}+\sum_{j\neq
i}l_j\frac{-b_i'b_j'}{b_ib_j}\right]$$
$$+(\overline{n}-1)[g(W,\nabla_VP)-\pi(V)\pi(W)].\eqno(3.5)$$ When $i\neq r$, then $\nabla_VP=\nabla_WP=\pi(V)=0$,
so
$$  \overline{{\rm Ric}} (V,W)={\rm Ric}^{F_i}
(V,W)+b_i^2g_{F_i}(V,W)\left[-\frac{b_i''}{b_i}+(l_i-1)\frac{-b_i'^2}{b_i^2}+\sum_{j\neq
i}l_j\frac{-b_i'b_j'}{b_ib_j}\right]=\lambda
b_i^2g_{F_i}(V,W).\eqno(3.6)$$ By variables separation, we have
$(F_i,\nabla^{F_i})~(i\neq r)$ is Einstein with the Einstein
constant $\lambda_i$ and
$$\lambda_i-b_ib_i''-b_ib_i'\sum_{j\neq
i}l_j\frac{b_j'}{b_j}-(l_i-1)b_i'^2=\lambda b_i^2.\eqno(3.7)$$ When
$i=r$, then $$ {\rm
Ric}^{F_i}(V,W)-g_{F_i}(V,W)\left[b_ib_i''+(l_i-1)b_i'^2+b_ib_i'\sum_{j\neq
i}l_j\frac{b_j'}{b_j}+\lambda b_i^2\right]$$
$$=(\overline{n}-1)\left[\pi(V)\pi(W)-\frac{g(W,\nabla_VP)+g(V,\nabla_WP)}{2}\right]\eqno(3.8)$$
So we prove the above theorem. ~~~~~$\Box$\\

When $M=I\times_{b_1}F_1\times_{b_2}F_2\cdots
    \times_{b_m}F_m$ be a multiply warped product and $P=\frac{\partial}{\partial
    t}$, by Proposition 2.10, we have
 $$
    \overline{S}=-2\sum_{i=1}^ml_i\frac{b_i''}{b_i}+\sum_{i=1}^m\frac{S^{F_i}}{b_i^2}+\sum_{i=1}^ml_i-\sum_{i=1}^ml_i(l_i-1)
\frac{b_i'^2}{b_i^2}-\sum_{i=1}^m\sum_{j\neq
i}l_il_j\frac{b_i'b_j'}{b_ib_j}+\sum_{i,j=1}^ml_il_j\frac{b_j'}{b_j} .\eqno(3.9)$$\\
The following result just follows from the method of separation of
variables and the fact that each $S^{F_i}$ is function defined on
$F_i$.\\

\noindent {\bf Proposition 3.4} {\it  Let
$M=I\times_{b_1}F_1\times_{b_2}F_2\cdots
    \times_{b_m}F_m$ be a multiply warped product and $P=\frac{\partial}{\partial
    t}$. If
    $(M,\overline{\nabla})$ has constant scalar curvature
    $\overline{S}$, then each $(F_i,\nabla^{F_i})$ has constant
    scalar curvature $S^{F_i}$.}\\

When $P\in\Gamma(TF_r)$, by Proposition 2.11, we have
$$
    \overline{S}=-2\sum_{i=1}^ml_i\frac{b_i''}{b_i}+\sum_{i=1}^m\frac{S^{F_i}}{b_i^2}+\sum_{i=1}^ml_i(l_i-1)
\frac{-b_i'^2}{b_i^2}+\sum_{i=1}^m\sum_{j\neq
i}l_il_j\frac{-b_i'b_j'}{b_ib_j}+
\pi(P)(1-\overline{n})+(\overline{n}-1){\rm div}_
{F_r}P.\eqno(3.10)$$\\

\noindent {\bf Proposition 3.5} {\it  Let
$M=I\times_{b_1}F_1\times_{b_2}F_2\cdots
    \times_{b_m}F_m$ be a multiply warped product and $P\in\Gamma(TF_r)$. If
    $(M,\overline{\nabla})$ has constant scalar curvature
    $\overline{S}$, then each $(F_i,\nabla^{F_i})~(i\neq r)$ has constant
    scalar curvature $S^{F_i}$ and if $g_{F_r}(P,P)$ and ${\rm div}_
{F_r}P$ are constants, then $S^{F_r}$ is also a constant. }\\

\section{ Generalized Robertson-Walker spacetimes with a
semi-symmetric non-metric connection}

\quad In this section, we study $M=I\times F$ with the metric tensor
$-dt^2+f(t)^2g_F$. As a corollary of Theorem 3.1, we
obtain:\\

\noindent {\bf Corollary 4.1} {\it Let $M=I\times F$ with the metric
tensor $-dt^2+f(t)^2g_F$ and $P=\frac{\partial}{\partial t}$. Then
    $(M,\overline{\nabla})$ is Einstein with the Einstein constant
    $\lambda$ if and only if the following conditions are satisfied
    }\\
\noindent (1){\it  $(F,\nabla^{F})$ is Einstein with the Einstein
constant $\lambda_F$.}\\
\noindent (2) $f''=(1-\frac{\lambda}{l})f.$\\ \noindent
(3)$\lambda_F+(1-l)f'^2+(\frac{\lambda}{l}-1-\lambda)f^2+lf'f=0.$\\

By Corollary 4.1 (2) and elementary methods for ordinary
differential
 equations, we get\\

\noindent {\bf Case i)} $\lambda<l$, then $f=c_1e^{at}+c_2e^{-at},$
where $a=\sqrt{1-\frac{\lambda}{l}}.$
 By Corollary
4.1 (3), then
$$\lambda_F+2c_1c_2\left[\frac{\lambda}{l}-1-\lambda+(l-1)a^2\right]+c_1^2e^{2at}\left[(1-l)a^2+(\frac{\lambda}{l}-1-\lambda)+la\right]$$
$$+c_2^2e^{-2at}\left[(1-l)a^2+(\frac{\lambda}{l}-1-\lambda)-la\right]=0.\eqno(4.1)$$
So
$$\lambda_F+2c_1c_2\left[\frac{\lambda}{l}-1-\lambda+(l-1)a^2\right]=0;\eqno(4.2i)$$
$$c_1^2\left[(1-l)a^2+(\frac{\lambda}{l}-1-\lambda)+la\right]=0;\eqno(4.2ii)$$
$$c_2^2\left[(1-l)a^2+(\frac{\lambda}{l}-1-\lambda)-la\right]=0.\eqno(4.2iii)$$
When $c_1\neq 0$, $c_2\neq 0$, by (4.2ii) and (4.2iii), then $la=0$,
this is a contradiction. When $c_1=0$ and $c_2\neq 0$, by (4.2iii),
$(1-l)a^2+(\frac{\lambda}{l}-1-\lambda)=la$. By
$a=\sqrt{1-\frac{\lambda}{l}}$, then $a=-1$, this is a
contradiction. When $c_1\neq 0$ and $c_2= 0$, by (4.2ii),
$(1-l)a^2+(\frac{\lambda}{l}-1-\lambda)=-la$, then $a=1$ and
$\lambda=0$. By (4.2i), $\lambda_F=0.$ So we get
$\underline{\lambda_F=\lambda=0,~f=c_1e^t.}$\\
\noindent {\bf Case ii)}$\lambda>l$, $f=c_1{\rm cos}(bt)+c_2{\rm
sin}(bt)$, where $b=\sqrt{\frac{\lambda}{l}-1}.$ By Corollary 4.1
(3), we have
$$\lambda_F+(1-l)(-c_1b{\rm sin}(bt)+c_2b{\rm cos}(bt))^2+(\frac{\lambda}{l}-1-\lambda)
(c_1{\rm cos}(bt)+c_2{\rm sin}(bt))^2$$ $$+l(c_1{\rm
cos}(bt)+c_2{\rm sin}(bt))(-c_1b{\rm sin}(bt)+c_2b{\rm
cos}(bt))=0.\eqno(4.3)$$ So
$$(1-l)c_1^2b^2+(\frac{\lambda}{l}-1-\lambda)c_2^2-lc_1c_2b=-\lambda_F;\eqno(4.3i)$$
$$(1-l)c_2^2b^2+(\frac{\lambda}{l}-1-\lambda)c_1^2+lc_1c_2b=-\lambda_F;\eqno(4.3ii)$$
$$-2(1-l)c_1c_2b^2+2c_1c_2(\frac{\lambda}{l}-1-\lambda)+l(-c_1^2+c_2^2)b=0.\eqno(4.3iii)$$
By (4.3i) and (4.3ii) and $b^2=\frac{\lambda}{l}-1$, we get
$$\left[(1-l)b^2-(\frac{\lambda}{l}-1-\lambda)\right](c_1^2-c_2^2)-2lc_1c_2b=0,~~c_1^2-c_2^2=2c_1c_2b.\eqno(4.4)$$
By (4.4) and (4.3iii), we have $c_1c_2=0$, then $c_1=c_2=0$. This is
a contradiction.\\
\noindent {\bf Case iii)} $\lambda=l$ and $f=c_1+c_2t$, by Corollary
4.1 (3), we get
$$\lambda_F+(1-l)c_2^2+(\frac{\lambda}{l}-1-\lambda)(c_1+c_2t)^2+l(c_1+c_2t)c_2=0.\eqno(4.5)$$
Then $$(\frac{\lambda}{l}-1-\lambda)c_2^2=0;\eqno(4.6i)$$
$$2(\frac{\lambda}{l}-1-\lambda)c_1c_2+lc_2^2=0;\eqno(4.6ii)$$
$$\lambda_F+(1-l)c_2^2+(\frac{\lambda}{l}-1-\lambda)c_1^2+lc_1c_2=0.\eqno(4.6iii)$$
By $\lambda=l$ and (4.6i), then $c_2=0$. By (4.6iii), then
$c_1=\sqrt{\frac{\lambda_F}{l}}.$ So we get
$\underline{\lambda=l,~f=\sqrt{\frac{\lambda_F}{l}}}.$ We get the
following theorem.\\

\noindent {\bf Theorem 4.2} {\it Let $M=I\times F$ with the metric
tensor $-dt^2+f(t)^2g_F$ and $P=\frac{\partial}{\partial t}$ and
${\rm dim}F>1$. Then
    $(M,\overline{\nabla})$ is Einstein with the Einstein constant
    $\lambda$ if and only if
 $(F,\nabla^{F})$ is Einstein with the Einstein
constant $\lambda_F$ and one of the following conditions holds
1)$\lambda_F=\lambda=0,~f=c_1e^t$ ~~2)$\lambda=l,~f=\sqrt{\frac{\lambda_F}{l}}$.}\\

\indent By (3.9) and Proposition 3.4, we have\\

\noindent {\bf Corollary 4.3} {\it  Let $M=I\times F$ with the
metric tensor $-dt^2+f(t)^2g_F$ and $P=\frac{\partial}{\partial t}$.
 If
    $(M,\overline{\nabla})$ has constant scalar curvature
    $\overline{S}$ if and only if $(F,\nabla^{F})$ has constant
    scalar curvature $S^{F}$ and}
$$\overline{S}=\frac{S^F}{f^2}-2l\frac{f''}{f}-l(l-1)\frac{f'^2}{f^2}+l+l^2\frac{f'}{f}.\eqno(4.7)$$\\

\indent In (4.7), we make the change of variable $f(t)=\sqrt{v(t)}$
and have the following equation
$$v''(t)+\frac{l-3}{4}\frac{v'(t)^2}{v(t)}-\frac{l}{2}v'(t)-\frac{l-\overline{S}}{l}v(t)-\frac{S^F}{l}=0.\eqno(4.8)$$\\

\noindent {\bf Theorem 4.4} {\it  Let $M=I\times F$ with the metric
tensor $-dt^2+f(t)^2g_F$ and $P=\frac{\partial}{\partial t}$ and
${\rm dim} F=l=3$.
 If
    $(M,\overline{\nabla})$ has constant scalar curvature
    $\overline{S}$ if and only if $(F,\nabla^{F})$ has constant
    scalar curvature $S^{F}$ and}\\
    \noindent (1) $\overline{S}<\frac{75}{16}$ and $\overline{S}\neq
    3$,~~$v(t)=c_1e^{\frac{\frac{3}{2}+\sqrt{\frac{25}{4}-\frac{4}{3}\overline{S}}}{2}t}+
    c_2e^{\frac{\frac{3}{2}-\sqrt{\frac{25}{4}-\frac{4}{3}\overline{S}}}{2}t}+\frac{S^F}{\overline{S}-3}.$\\
\noindent (2) $\overline{S}=\frac{75}{16},~~
~~v(t)=c_1e^{\frac{3}{4}t}+c_2te^{\frac{3}{4}t}+\frac{S^F}{\overline{S}-3}.$\\
\noindent (3) $\overline{S}>\frac{75}{16},~~
~~v(t)=c_1e^{\frac{3}{4}t}{\rm
cos}\left(\frac{\sqrt{\frac{3}{4}\overline{S}-\frac{25}{4}}}{2}t\right)
+c_2e^{\frac{3}{4}t}{\rm
sin}\left(\frac{\sqrt{\frac{3}{4}\overline{S}-\frac{25}{4}}}{2}t\right)+\frac{S^F}{\overline{S}-3}.$\\
\noindent (4)
$\overline{S}=3,~~~~v(t)=c_1-\frac{2S^F}{9}t+c_2e^{\frac{3}{2}t}.$\\

\noindent{\bf Proof.} If $l=3$, then we have a simple differential
equation
$$v''(t)-\frac{3}{2}v'(t)+(\frac{\overline{S}}{3}-1)v(t)-\frac{S^F}{3}=0.\eqno(4.9)$$\\
If $\overline{S}\neq
    3$, we putting
    $h(t)=(\frac{\overline{S}}{3}-1)v(t)-\frac{S^F}{3},$ it follows
    that $h''(t)-\frac{3}{2}h'(t)+(\frac{\overline{S}}{3}-1)h(t)=0$. The above
    solutions (1)-(3) follow directly from elementary methods for
    ordinary differential equations. When $\overline{S}=
    3$, then $v''(t)-\frac{3}{2}v'(t)-\frac{S^F}{3}=0$, we get the solution
    (4).$\Box$\\

\noindent {\bf Theorem 4.5} {\it  Let $M=I\times F$ with the metric
tensor $-dt^2+f(t)^2g_F$ and $P=\frac{\partial}{\partial t}$ and
${\rm dim} F=l\neq 3$ and $S^F=0$. Let
$\triangle=\frac{l^2}{4}+\frac{(l+1)(l-\overline{S})}{l}.$
 If
    $(M,\overline{\nabla})$ has constant scalar curvature
    $\overline{S}$ if and only if }\\
    \noindent (1)
    $\overline{S}<\frac{l^3}{4(l+1)}+l,~~v(t)=\left(c_1e^{\frac{\frac{l}{2}+\sqrt{\triangle}}{2}t}
    +c_2e^{\frac{\frac{l}{2}-\sqrt{\triangle}}{2}t}\right)^{\frac{4}{l+1}}.$\\

 \noindent (2)$\overline{S}=\frac{l^3}{4(l+1)}+l
   ,~~v(t)=\left(c_1e^{\frac{l}{4}t}
    +c_2te^{\frac{l}{4}t}\right)^{\frac{4}{l+1}}.$\\

\noindent (3)$\overline{S}>\frac{l^3}{4(l+1)}+l
   ,~~v(t)=\left(c_1e^{\frac{l}{4}t}{\rm
    cos}\left(\frac{\sqrt{-\triangle}}{2}t\right)
    +c_2e^{\frac{l}{4}t}{\rm
    sin}\left(\frac{\sqrt{-\triangle}}{2}t\right)\right)^{\frac{4}{l+1}}.$\\

\noindent{\bf Proof.} In this case, the equation (4.8) is changed
into the simpler form
$$\frac{v''(t)}{v(t)}+\frac{l-3}{4}\frac{v'(t)^2}{v(t)^2}-\frac{l}{2}\frac{v'(t)}{v(t)}-\frac{l-\overline{S}}{l}=0.\eqno(4.10)$$\\
Putting $v(t)=w(t)^{\frac{4}{l+1}}$, then $w(t)$ satisfies the
equation
$w''-\frac{l}{2}w'+\frac{(l+1)}{4}\frac{(\overline{S}-l)}{l}w=0$, by
the elementary methods for
    ordinary differential equations, we prove the above
    theorem.~~~~~~~~~$\Box$\\

When ${\rm dim} F=l\neq 3$ and $S^F\neq 0$, putting
$v(t)=w(t)^{\frac{4}{l+1}}$, then $w(t)$ satisfies the equation
$$w''-\frac{l}{2}w'+\frac{(l+1)}{4}\frac{(\overline{S}-l)}{l}w-\frac{(l+1)}{4}\frac{S^F}{l}w^{1-\frac{4}{l+1}}=0.\eqno(4.11)$$

\section{Generalized Kasner spacetimes with a
semi-symmetric non-metric
 connection }

 \quad In this section, we consider the scalar and Ricci
 curvature of generalized Kasner spacetimes with a semi-symmetric non-metric
 connection. We recall the definition of generalized Kasner
 spacetimes ([DU]).\\

 \noindent{\bf Definition 5.1} A generalized Kasner spacetime
 $(M,g)$ is a Lorentzian multiply warped product of the form
 $M=I\times_{\phi^{p_1}}F_1\times\cdots\times_{\phi^{p_m}}F_m$ with
 the metric
 $g=-dt^2\oplus\phi^{2p_1}g_{F_1}\oplus\cdots\oplus\phi^{2p_m}g_{F_m}$,
 where $\phi:~I\rightarrow (0,\infty)$ is smooth and
 $p_i\in\mathbb{R}$, for any $i\in\{1,\cdots,m\}$ and also
 $I=(t_1,t_2)$.\\

\indent We introduce the following parameters
$\zeta=\sum_{i=1}^ml_ip_i$ and $\eta=\sum_{i=1}^ml_ip_i^2$ for
generalized Kasner spacetimes. By Theorem 3.1 and direct
computations, we get\\

\noindent {\bf Proposition 5.2}{\it~Let
$M=I\times_{\phi^{p_1}}F_1\times\cdots\times_{\phi^{p_m}}F_m$ be a
generalized Kasner spacetime and $P=\frac{\partial}{\partial t}$.
Then
    $(M,\overline{\nabla})$ is Einstein with the Einstein constant
    $\lambda$ if and only if the following conditions are satisfied
    for any $i\in\{1,\cdots,m\}$}\\
\noindent (1){\it  $(F_i,\nabla^{F_i})$ is Einstein with the
Einstein
constant $\lambda_i$, $i\in\{1,\cdots,m\}$.}\\
\noindent (2)
$(\eta-\zeta)\frac{\phi'^2}{\phi^2}+\zeta\left(\frac{\phi''}{\phi}\right)+\lambda-\sum_{i=1}^ml_i=0.$\\
\noindent
(3)$\frac{\lambda_i}{\phi^{2p_i}}-p_i\frac{\phi''}{\phi}-(\zeta-1)p_i\frac{\phi'^2}{\phi^2}+\zeta\frac{\phi'}{\phi}=
\lambda.$\\

By (3.9) we obtain\\

\noindent {\bf Proposition 5.3}{\it~Let
$M=I\times_{\phi^{p_1}}F_1\times\cdots\times_{\phi^{p_m}}F_m$ be a
generalized Kasner spacetime and $P=\frac{\partial}{\partial t}$.
Then
    $(M,\overline{\nabla})$ has constant scalar curvature
    $\overline{S}$ if and only if each $(F_i,\nabla^{F_i})$ has constant
    scalar curvature $S^{F_i}$ and}
$$\overline{S}=\sum_{i=1}^m\frac{S^{F_i}}{\phi^{2p_i}}-2\zeta\frac{\phi''}{\phi}-(\eta+\zeta^2-2\zeta)\frac{\phi'^2}{\phi^2}+(\overline{n}-1)\zeta
\frac{\phi'}{\phi}+(\overline{n}-1).\eqno(5.1)$$

\indent Nextly, we first give a classification of four-dimensional
generalized Kasner spacetimes with a semi-symmetric non-metric
 connection and then consider Ricci tensors and scalar curvatures
 of them.\\

\noindent{\bf Definition 5.4} Let
 $M=I\times_{b_1}F_1\times\cdots\times_{b_m}F_m$ with
 the metric
 $g=-dt^2\oplus b_1^2g_{F_1}\oplus\cdots\oplus b_m^2g_{F_m}$.\\
\noindent {\bf $\cdot$} $(M,g)$ is said to be of Type (I) if $m=1$
and ${\rm dim}(F)=3$.\\
\noindent {\bf $\cdot$} $(M,g)$ is said to be of Type (II) if $m=2$
and ${\rm dim}(F_1)=1$ and ${\rm dim}(F_2)=2$.\\
\noindent {\bf $\cdot$} $(M,g)$ is said to be of Type (III) if $m=3$
and ${\rm dim}(F_1)=1$,~${\rm dim}(F_2)=1$ and ${\rm dim}(F_3)=1$.\\

By Theorem 4.2 and 4.4, we have given a classification of Type (I)
Einstein spaces and Type (I) spaces with the constant scalar
curvature.\\

\noindent{\bf $\cdot$ Classification of Einstein Type (II)
generalized Kasner space-times with a semi-symmetric non-metric
 connection}\\
 \indent Let $M=I\times_{\phi^{p_1}}F_1\times_{\phi^{p_2}}F_2$ be an
 Einstein type (II) generalized Kasner spacetime and $P=\frac{\partial}{\partial
 t}$. Then $\zeta=p_1+2p_2$, $\eta=p_1^2+2p_2^2$.  By Proposition
 5.2
 , we have
$$(\eta-\zeta)\frac{\phi'^2}{\phi^2}+\zeta\left(\frac{\phi''}{\phi}\right)+\lambda-3=0;\eqno(5.2i)$$
$$-p_1\frac{\phi''}{\phi}-(\zeta-1)p_1\frac{\phi'^2}{\phi^2}+\zeta\frac{\phi'}{\phi}=
\lambda;\eqno(5.2ii)$$
$$\frac{\lambda_2}{\phi^{2p_2}}-p_2\frac{\phi''}{\phi}-(\zeta-1)p_2\frac{\phi'^2}{\phi^2}+\zeta\frac{\phi'}{\phi}=
\lambda.\eqno(5.2iii)$$
 where $\lambda_2$ is a constant. Consider
following two cases:\\

\noindent {\bf Case i)}~ $\underline{\zeta=0}$\\
 \indent In this
case, $p_2=-\frac{1}{2}p_1$, $\eta=\frac{3}{2}p_1^2$. Then by (5.2),
we have
$$\eta\frac{\phi'^2}{\phi^2}+\lambda-3=0 ,\eqno(5.3i)$$

        $$ p_1\left(-\frac{\phi''}{\phi}+\frac{\phi'^2}{\phi^2}\right)=
\lambda,\eqno(5.3ii)$$
       $$ \frac{\lambda_2}{\phi^{-p_1}}-\frac{1}{2}p_1\left(-\frac{\phi''}{\phi}+\frac{\phi'^2}{\phi^2}\right)
      =
\lambda,\eqno(5.3iii)$$
 {\bf Case i a)}~ $\underline{\eta=0}$, then $p_i=0$, by (5.3i), $\lambda=3$. By (5.3ii), $\lambda=0$, this
is a contradiction.
\\
{\bf Case i b)}~ $\underline{\eta\neq 0}$, then $p_i\neq 0$.\\
 {\bf Case i b)1)} $\underline{\lambda_2=0}$\\
by (5.3ii) and (5.3iii), $\lambda=0$ and
$$-\frac{\phi''}{\phi}+\frac{\phi'^2}{\phi^2}=0,~~
\frac{\phi'^2}{\phi^2}=\frac{3}{\eta},\eqno(5.4)$$ then
$\phi=c_0e^{\pm\sqrt{\frac{3}{\eta}}t}$ which satisfies the first
equation in (5.4), so we get\\ $\underline{\zeta=0, ~\eta\neq
0,~\lambda=\lambda_F=0,~\phi=c_0e^{\pm\sqrt{\frac{3}{\eta}}t}}.$
\\
{\bf Case i b)2)} $\underline{\lambda_2\neq 0}$\\
by (5.3ii) and (5.3iii), we have
$\frac{\lambda_2}{\phi^{-p_1}}=\frac{3}{2}\lambda,$ so $\phi$ is a
constant. By (5.3ii), $\lambda=0$, so $\lambda_2=0$, this is a
contradiction.\\

\noindent {\bf Case ii)}~ $\underline{\zeta\neq  0},$\\
then $\eta\neq 0$. Putting $\phi=\psi^{\frac{\zeta}{\eta}}$, then
$\psi''+\frac{(\lambda-3)\eta}{\zeta^2}\psi=0$. Hence,\\

\noindent (1)~~$\lambda<3,$
$\psi=c_1e^{\sqrt{\frac{(3-\lambda)\eta}{\zeta^2}}t} +
c_2e^{-\sqrt{\frac{(3-\lambda)\eta}{\zeta^2}}t},$\\
\noindent (2)~~$\lambda=3,$
$\psi=c_1 + c_2t,$\\
\noindent (3)~~$\lambda>3,$ $\psi=c_1{\rm
cos}\left(\sqrt{\frac{(\lambda-3)\eta}{\zeta^2}}t \right)
 + c_2{\rm sin}\left(\sqrt{\frac{(\lambda-3)\eta}{\zeta^2}}t
\right).$\\
We make (5.2) into
$$\psi''+\frac{(\lambda-3)\eta}{\zeta^2}\psi=0;\eqno(5.5i)$$
$$-\frac{p_1}{\zeta}\frac{(\phi^\zeta)''}{\phi^\zeta}+
\frac{(\phi^\zeta)'}{\phi^\zeta}=\lambda; \eqno(5.5ii)$$
$$\frac{\lambda_2}{\phi^{2p_2}}-\frac{p_2}{\zeta}\frac{(\phi^\zeta)''}{\phi^\zeta}+
\frac{(\phi^\zeta)'}{\phi^\zeta}=\lambda.\eqno(5.5iii)$$ When
$p_1=p_2$, the type (II) spaces turns into type (I) spaces, so we
assume $p_1\neq p_2$.
 By (5.5ii) and (5.5iii), then
$$\psi'=\frac{p_1\lambda_2\eta}{(p_2-p_1)\zeta^2}\psi^{1-\frac{2p_2\zeta}{\eta}}
+ \frac{\lambda\eta}{\zeta^2}\psi.\eqno(5.6)$$
 \noindent {\bf Case ii)(1)}~~$\lambda<3,$
 $\psi=c_1e^{at} +
c_2e^{-at},$ where $a=\sqrt{\frac{(3-\lambda)\eta}{\zeta^2}}.$\\
By (5.6),
$$ac_1e^{at} - ac_2e^{-at}=\frac{p_1\lambda_2\eta}{(p_2-p_1)\zeta^2}(c_1e^{at} + c_2e^{-at})^{1-\frac{2p_2\zeta}{\eta}}
+ \frac{\lambda\eta}{\zeta^2}(c_1e^{at} + c_2e^{-at}).\eqno(5.7)$$
 \noindent {\bf Case ii)(1)(a)}~$\underline{c_1=0}$,\\
 then $$-\left[a+\frac{\lambda\eta}{\zeta^2}\right]c_2e^{-at}=
\frac{p_1\lambda_2\eta}{(p_2-p_1)\zeta^2}(
c_2e^{-at})^{1-\frac{2p_2\zeta}{\eta}}.\eqno(5.8)$$
 \noindent {\bf Case ii)(1)(a)1)}~$~\underline{p_1\lambda_2\neq
 0}$\\
 then $p_2=0$ and $\zeta=p_1,~\eta=p_1^2$ and
  $\psi=c_2e^{-at}$. By
 (5.5ii), we get $-a^2-a=\lambda$ and
 $-\sqrt{3-\lambda}=3$, this is a contradiction.\\
\noindent {\bf Case ii)(1)(a)2)}~$\underline{p_1\lambda_2=
 0}$\\
If $\underline{p_1=0},$ then $\zeta=2p_2$, $\eta=2p_2^2$. By (5.8),
$a=-\frac{\lambda}{2}$ and $\phi^\zeta=c_2e^{-2at}.$ By (5.5iii),
then
$$\frac{\lambda_2}{c'e^{\frac{-
4p_2at}{\zeta}}}-2a^2-2a=\lambda,\eqno(5.9)$$ so $\lambda_2=0$ and
$-2a^2-2a=\lambda$. By $a=-\frac{\lambda}{2}$, then $\lambda=0$ and
$a=0$, this is a contradiction.\\
If $\underline{\lambda_2=0,~p_1\neq 0}$, by (5.5ii) and
$a=-\frac{\lambda\eta}{\zeta^2}$ and
$\phi^\zeta=c'e^{-\frac{at\zeta^2}{\eta}},$ then $\lambda=0$ and
$a=0$, this is a
contradiction.\\

 \noindent {\bf Case ii)(1)(b)}~$\underline{c_2=0}$,\\
then $$\left[a-\frac{\lambda\eta}{\zeta^2}\right]c_1e^{at}=
\frac{p_1\lambda_2\eta}{(p_2-p_1)\zeta^2}(
c_1e^{at})^{1-\frac{2p_2\zeta}{\eta}}.\eqno(5.10)$$

 \noindent {\bf
Case ii)(1)(b)1)} $\underline{p_1\lambda_2}\neq 0$,\\
then $p_2=0$ and $\zeta=p_1$, $\eta=p_1^2$ and
$\sqrt{3-\lambda}-\lambda=-\lambda_2.$ By (5.5ii), we get
$-a^2+a=\lambda$ and $\lambda=-6$, so $\lambda_2=-9$. In this case,
(5.2iii) holds, so we get $\underline{p_2=0,~p_1\neq
0,~\lambda=-6,~\lambda_2=-9,~\phi=c_1e^{\frac{3t}{\zeta}}}.$

\noindent {\bf
Case ii)(1)(b)2)} $\underline{p_1\lambda_2}= 0$,\\
if $\underline{p_1=0}$, then $\zeta=2p_2$, $\eta=2p_2^2$ and
$\psi=c_1e^{at}$ and
$a=\frac{\lambda\eta}{\zeta^2}=\frac{\lambda}{2}$. By (5.5iii),
$$\frac{\lambda_2}{c'e^{\frac{
4p_2at}{\zeta}}}-2a^2+2a=\lambda,\eqno(5.11)$$ so $\lambda_2=0$ and
$-2a^2+2a=\lambda$, then $\lambda=a=0$, this is a contradiction. If
$\underline{\lambda_2=0}$, by (5.5ii), then $\lambda=a=0$, this
is a contradiction.\\

\noindent {\bf Case ii)(1)(c)}~$\underline{c_1\neq 0,~c_2\neq
0},$\\
If $\underline{p_2\neq 0}$, then
$e^{at},e^{-at},(c_1e^{at}+c_2e^{-at})^{1-\frac{2p_2\zeta}{\eta}}$
are linear independent, by (5.7), then
$$\left[a-\frac{\lambda\eta}{\zeta^2}\right]c_1=0,~\left[-a-\frac{\lambda\eta}{\zeta^2}\right]c_2=0,~
\frac{p_1\lambda_2\eta}{(p_2-p_1)\zeta^2}=0.\eqno(5.12)$$ So
$a=0$, this is a contradiction.\\
If $\underline{p_2= 0}$, then by (5.7),
$$a-\frac{\lambda\eta}{\zeta^2}-\frac{p_1\lambda_2\eta}{(p_2-p_1)\zeta^2}=0,~
-a-\frac{\lambda\eta}{\zeta^2}-\frac{p_1\lambda_2\eta}{(p_2-p_1)\zeta^2}=0,\eqno(5.13)$$
so $a=0$ and we get a contradiction.\\

\noindent {\bf Case ii)(2)}~~$\lambda=3,$
$\psi=c_1+ c_2t,$\\
by (5.6), we have
$$c_2=\frac{p_1\lambda_2\eta}{(p_2-p_1)\zeta^2}(c_1+c_2t)^{1-\frac{2p_2\zeta}{\eta}}
+ \frac{3\eta}{\zeta^2}(c_1+c_2t).\eqno(5.14)$$
 \noindent {\bf Case ii)(2)a)}~~$\underline{c_2\neq 0},$\\
so $p_2=0$. By (5.5iii), then $
\phi^\zeta=c_0e^{(-\lambda_2+\lambda)t}$ and $c_1 +
c_2t=c_0e^{(-\lambda_2+\lambda)t}$,
this is a contradiction with $c_2\neq 0.$\\
\noindent {\bf Case ii)(2)b)}~~$\underline{c_2= 0},$\\
then $\psi$ and $\phi$ are constants, by (5.5ii), then $\lambda=0$.
This is a contradiction with $\lambda=3.$\\

\noindent {\bf Case ii)(3)}~~$\lambda>3,$ $\psi=c_1{\rm cos}\left(a
t\right)
 + c_2{\rm sin}\left(a
t\right),$ where $a=\sqrt{\frac{(\lambda-3)\eta}{\zeta^2}}$. By
(5.6), we have
$$a(-c_1{\rm sin}(at)+c_2{\rm
cos}(at))=\frac{p_1\lambda_2\eta}{(p_2-p_1)\zeta^2}(c_1{\rm
cos}(at)+c_2{\rm sin}(at))^{1-\frac{2p_2\zeta}{\eta}}$$ $$ +
\frac{\lambda\eta}{\zeta^2}(c_1{\rm cos}(at)+c_2{\rm
sin}(at)).\eqno(5.15)$$
 If $\underline{p_2\neq 0},$ then ${\rm sin}\left(a
t\right)$, ${\rm cos}\left(a t\right)$,~$(c_1{\rm cos}(at)+c_2{\rm
sin}(at))^{1-\frac{2p_2\zeta}{\eta}}$ are linear indepdent, so
$-ac_1=\frac{\lambda\eta}{\zeta^2}c_2,~ac_2=\frac{\lambda\eta}{\zeta^2}c_1$
and $a=0$, this is a contradiction. If $\underline{p_2= 0},$ then
$$-ac_1=\frac{p_1\lambda_2\eta}{(p_2-p_1)\zeta^2}c_2+\frac{\lambda\eta}{\zeta^2}c_2;~~
ac_2=\frac{p_1\lambda_2\eta}{(p_2-p_1)\zeta^2}c_1+\frac{\lambda\eta}{\zeta^2}c_1.\eqno(5.16)$$
Then $a=0$, this is a contradiction. By the above discussions, we get the following theorem:\\

\noindent {\bf Theorem 5.5}{\it~Let
$M=I\times_{\phi^{p_1}}F_1\times_{\phi^{p_2}}F_2$ be a generalized
Kasner spacetime and ${\rm dim}F_1=1,~{\rm dim}F_2=2$ and
$P=\frac{\partial}{\partial t}$. Then
    $(M,\overline{\nabla})$ is Einstein with the Einstein constant
    $\lambda$ if and only if
 $(F_2,\nabla^{F_2})$ is Einstein with the Einstein
constant $\lambda_2$, and  one of the following conditions is satisfied}\\
\noindent (1)$\zeta=0, ~\eta\neq
0,~\lambda=\lambda_F=0,~\phi=c_0e^{\pm\sqrt{\frac{3}{\eta}}t}.$\\
\noindent (2)$p_2=0,~p_1\neq
0,~\lambda=-6,~\lambda_2=-9,~\phi=c_1e^{\frac{3t}{\zeta}}.$\\

\noindent{\bf $\cdot$ Type (II) generalized Kasner space-times with
a semi-symmetric non-metric
 connection with constant scalar curvature}\\
\indent By Proposition 5.3, then $(F_2,\nabla^{F_2})$ has constant
    scalar curvature $S^{F_2}$ and
$$\overline{S}=\frac{S^{F_2}}{\phi^{2p_2}}-2\zeta\frac{\phi''}{\phi}-(\eta+\zeta^2-2\zeta)\frac{\phi'^2}{\phi^2}+3\zeta
\frac{\phi'}{\phi}+3.\eqno(5.17)$$
 If $\underline{\zeta=0}$, when
$\underline{\eta=0},$  then $p_1=p_2=0$ and
$\underline{\overline{S}=S^{F_2}+3}$. If $\underline{\eta\neq 0},$
then
$$\eta\frac{\phi'^2}{\phi^2}=\frac{S^{F_2}}{\phi^{2p_2}}+(-\overline{S}+3).\eqno(5.18)$$
If $\underline{\zeta\neq 0},$ putting
$\phi=\psi^{\frac{2\zeta}{\eta+\zeta^2}},$ we get
$$-\frac{4\zeta^2}{\eta+\zeta^2}\psi''+\frac{6\zeta^2}{\eta+\zeta^2}\psi'+(-\overline{S}+3)\psi
+S^{F_2}\psi^{1-\frac{4p_2\zeta}{\eta+\zeta^2}}=0.\eqno(5.19)$$

\noindent{\bf $\cdot$ Type (III) generalized Kasner space-times with
a semi-symmetric non-metric
 connection with constant scalar curvature}\\

\indent By Proposition 5.3, then
$$\overline{S}=-2\zeta\frac{\phi''}{\phi}-(\eta+\zeta^2-2\zeta)\frac{\phi'^2}{\phi^2}+3\zeta
\frac{\phi'}{\phi}+3.\eqno(5.20)$$ If $\underline{\zeta=\eta=0},$
then $p_1=p_2=p_3=0$, we get $\overline{S}=3$.\\
If $\underline{\zeta=0,~\eta\neq 0},$ then $[({\rm
ln}\phi)']^2=-\frac{\overline{S}-3}{\eta}$, so when
$\overline{S}>3$, there is no solutions, when $\overline{S}=3$,
$\phi$ is a constant and  when $\overline{S}<3$,
$\phi=c_0e^{\pm\sqrt{\frac{-\overline{S}+3}{\eta}}t}.$\\
If $\underline{\zeta\neq 0}$, then $\eta\neq 0$, putting
$\phi=\psi^{\frac{2\zeta}{\eta+\zeta^2}}$, then
$$\psi''-\frac{3}{2}\psi'+\frac{(\overline{S}-3)(\eta+\zeta^2)}{4\zeta^2}\psi=0.\eqno(5.21)$$
Let
$\triangle=\frac{9}{4}-\frac{(\overline{S}-3)(\eta+\zeta^2)}{\zeta^2}.$
So we get \\

\noindent (1) $\overline{S}<\frac{9\zeta^2}{4(\eta+\zeta^2)}+3,$~
$\psi=c_1e^{\frac{\frac{3}{2}+\sqrt{\triangle}}{2}t}+c_2
e^{\frac{\frac{3}{2}-\sqrt{\triangle}}{2}t},$\\

\noindent (2)$\overline{S}=\frac{9\zeta^2}{4(\eta+\zeta^2)}+3,$~
$\psi=c_1e^{\frac{3}{4}t}+c_2te^{\frac{3}{4}t},$\\
\noindent (3) $\overline{S}>\frac{9\zeta^2}{4(\eta+\zeta^2)}+3,$~
$\psi=c_1e^{\frac{3}{4}t}{\rm
cos}\left(\frac{\sqrt{-\triangle}}{2}t\right)+c_2e^{\frac{3}{4}t}
{\rm sin}\left(\frac{\sqrt{-\triangle}}{2}t\right).$
So we get the following theorem\\

\noindent {\bf Theorem 5.6}{\it~Let
$M=I\times_{\phi^{p_1}}F_1\times_{\phi^{p_2}}F_2\times_{\phi^{p_3}}F_3$
be a generalized Kasner spacetime and ${\rm dim}F_1={\rm
dim}F_2={\rm dim}F_3=1,$ and $P=\frac{\partial}{\partial t}$. Then
$\overline{S}$ is a constant if and only if one of the following
case holds} \\
\noindent (1) $\zeta=\eta=0,$ $\overline{S}=3$.\\
\noindent (2){\it $ \zeta=0,~\eta\neq 0,$ when $\overline{S}> 3$,
there is no solutions, when $\overline{S}=3$, $\phi$ is a constant
and
when $\overline{S}<3$, $\phi=c_0e^{\pm\sqrt{-\frac{\overline{S}-3}{\eta}}t}.$}\\
\noindent (3) If $\zeta\neq 0$\\
 \noindent (3a) $\overline{S}<\frac{9\zeta^2}{4(\eta+\zeta^2)}+3,$~
$\phi=\left(c_1e^{\frac{\frac{3}{2}+\sqrt{\triangle}}{2}t}+c_2
e^{\frac{\frac{3}{2}-\sqrt{\triangle}}{2}t}\right)^{\frac{2\zeta}{\eta+\zeta^2}},$\\
\noindent (3b) $\overline{S}=\frac{9\zeta^2}{4(\eta+\zeta^2)}+3,$~
$\phi=\left(c_1e^{\frac{3}{4}t}+c_2te^{\frac{3}{4}t}\right)^{\frac{2\zeta}{\eta+\zeta^2}},$\\
\noindent (3c) $\overline{S}>\frac{9\zeta^2}{4(\eta+\zeta^2)}+3,$~
$\phi=\left(c_1e^{\frac{3}{4}t}{\rm
cos}\left(\frac{\sqrt{-\triangle}}{2}t\right)+c_2e^{\frac{3}{4}t}
{\rm sin}\left(\frac{\sqrt{-\triangle}}{2}t\right)
\right)^{\frac{2\zeta}{\eta+\zeta^2}}.$\\

\noindent{\bf $\cdot$ Einstein Type (III) generalized Kasner
space-times with a semi-symmetric non-metric
 connection}\\

\indent By Proposition 5.2, we have\\
$$\zeta\left(\frac{\phi''}{\phi}\right)+(\eta-\zeta)\frac{\phi'^2}{\phi^2}+\lambda-3=0,\eqno(5.22i)$$
$$-p_1\left[\frac{\phi''}{\phi}+(\zeta-1)\frac{\phi'^2}{\phi^2}\right]+\zeta\frac{\phi'}{\phi}=\lambda,\eqno(5.22ii)$$
$$-p_2\left[\frac{\phi''}{\phi}+(\zeta-1)\frac{\phi'^2}{\phi^2}\right]+\zeta\frac{\phi'}{\phi}=\lambda,\eqno(5.22iii)$$
$$-p_3\left[\frac{\phi''}{\phi}+(\zeta-1)\frac{\phi'^2}{\phi^2}\right]+\zeta\frac{\phi'}{\phi}=\lambda,\eqno(5.22iv)$$
If $\underline{\zeta=\eta=0},$ by (5.22i), $\lambda=3$, by (4.50ii),
$\lambda=0$, this is a contradiction.\\
If $\underline{\zeta=0, ~\eta\neq0},$ plusing
(5.22ii),(5.22iii),(5.22iv), we get $\lambda=0$. By (5.22i),
$\frac{\phi'^2}{\phi^2}=\frac{3}{\eta}$ and
$\phi=c_0e^{\pm\sqrt{\frac{3}{\eta}}t}$ which satisfies (5.22ii),
(5.22iii) and (5.22iv), so we obtain
$\underline{\lambda=0,~\zeta=0,~\eta\neq
0,~\phi=c_0e^{\pm\sqrt{\frac{3}{\eta}}t}}.$\\
When $\underline{\zeta\neq 0}$, if $p_1=p_2=p_3$, we get type (I),
so we may let $p_1\neq p_2$. By (5.22ii) and (5.22iii), we have
$\frac{(\phi^\zeta)'}{\phi^\zeta}=\lambda$ and
$\frac{(\phi^\zeta)''}{\phi^\zeta}=0,$ so $\phi^\zeta=c_1+c_2t$ and
$c_2=\lambda(c_1+c_2t),$ then $\lambda c_2=0$. When $c_2=0$, then
$\lambda=0$ and $\phi$ is a constant, so by (5.22i), then
$\lambda=3$ which is a contradiction. When $\lambda=0,$ then $c_2=0$
which is also a contradiction. We get the following
theorem.\\

\noindent {\bf Theorem 5.7}{\it~Let
$M=I\times_{\phi^{p_1}}F_1\times_{\phi^{p_2}}F_2\times_{\phi^{p_3}}F_3$
be a generalized Kasner spacetime for $p_i\neq p_j$ for some
$i,j\in\{1,2,3\}$ and ${\rm dim}F_1={\rm dim}F_2={\rm dim}F_3=1,$
and $P=\frac{\partial}{\partial t}$.
 Then
    $(M,\overline{\nabla})$ is Einstein with the Einstein constant
    $\lambda$ if and only if $\lambda=0,~\zeta=0,~\eta\neq
0,~\phi=c_0e^{\pm\sqrt{\frac{3}{\eta}}t}.$}\\

\section{ Curvature
      of multiply twisted products with an affine
 connection with a zero torsion}

\quad Let $\nabla$ is the Levi-Civita connection of $M$, we define
$$\widetilde{\nabla}_XY=\nabla_XY+\pi(X)Y+\pi(Y)X,\eqno(6.1)$$
which has no torsion. Let $R$ and $\widetilde{R}$ be the curvature
tensors of $\nabla$ and $\widetilde{\nabla}$ respectively. Then $R$
and $\widetilde{R}$ are related by
$$\widetilde{R}(X,Y)Z=R(X,Y)Z+X(\pi(Y))Z-Y(\pi(X))Z$$
$$
+g(Z,\nabla_XP)Y-g(Z,\nabla_YP)X+
\pi(Z)[\pi(Y)X-\pi(X)Y]-\pi([X,Y])Z$$
$$=\overline{R}(X,Y)Z+X(\pi(Y))Z-Y(\pi(X))Z-\pi([X,Y])Z,\eqno(6.2)$$ for any vector fields $X,Y,Z$ on $M$. By (6.1) and Proposition
2.2 in [Wa], we have\\

 \noindent {\bf Proposition 6.1} {\it Let $M=B\times_{b_1}F_1\times_{b_2}F_2\cdots
    \times_{b_m}F_m$ be a multiply twisted product and let $X,Y\in \Gamma(TB)$
      and $U\in \Gamma(TF_i)$, $W\in \Gamma(TF_j)$ and $P\in \Gamma(TB)$ .
      Then}
      \noindent$(1) ~~\widetilde{\nabla}_XY=\widetilde{\nabla}^B_XY.$\\
      \noindent$(2)~~\widetilde{\nabla}_XU=[\frac{X(b_i)}{b_i}+\pi(X)]U.$\\
\noindent$(3)~~\widetilde{\nabla}_UX=[\frac{X(b_i)}{b_i}+\pi(X)]U.$\\
            \noindent$(4)~~\widetilde{\nabla}_UW=0~~if~ i\neq j.$\\
      \noindent$(5)~~\widetilde{\nabla}_UW=U({\rm ln}b_i)W+W({\rm ln}b_i)U-\frac{g_{F_i}(U,W)}{b_i}{\rm
      grad}_{F_i}b_i-b_ig_{F_i}(U,W){\rm
      grad}_{B}b_i+\nabla^{F_i}_UW~~if~ i= j.$\\

\noindent {\bf Proposition 6.2} {\it Let
$M=B\times_{b_1}F_1\times_{b_2}F_2\cdots
    \times_{b_m}F_m$ be a multiply twisted product and let $X,Y\in \Gamma(TB)$
      and $U\in \Gamma(TF_i)$, $W\in \Gamma(TF_j)$ and $P\in \Gamma(TF_k)$ . Then}\\
      \noindent$(1) ~~\widetilde{\nabla}_XY={\nabla}^B_XY.$\\
      \noindent$(2)~~\widetilde{\nabla}_XU=\frac{X(b_i)}{b_i}U+g(P,U)X.$\\
\noindent$(3)~~\widetilde{\nabla}_UX=\frac{X(b_i)}{b_i}U+g(P,U)X.$\\
            \noindent$(4)~~\widetilde{\nabla}_UW=\pi(W)U+\pi(U)W~~if~ i\neq j.$\\
      \noindent$(5)~~\widetilde{\nabla}_UW=U({\rm ln}b_i)W+W({\rm ln}b_i)U-\frac{g_{F_i}(U,W)}{b_i}{\rm
      grad}_{F_i}b_i-b_ig_{F_i}(U,W){\rm
      grad}_{B}b_i+\nabla^{F_i}_UW+\pi(W)U+\pi(U)W~~if~ i= j.$\\

\indent By (6.2) and Proposition 2.4 and 2.5, we have\\

\noindent {\bf Proposition 6.3} {\it Let
$M=B\times_{b_1}F_1\times_{b_2}F_2\cdots
    \times_{b_m}F_m$ be a multiply twisted product and let $X,Y,Z\in \Gamma(TB)$
      and $V\in \Gamma(TF_i)$, $W\in \Gamma(TF_j)$, $U\in \Gamma(TF_k)$ and $P\in \Gamma(TB)$. Then}\\
\noindent $(1)\widetilde{R}(X,Y)Z=\widetilde{R}^B(X,Y)Z,$\\
\noindent $(2)\widetilde{R}(X,Y)V=X(\pi(Y))V-Y(\pi(X))V-\pi([X,Y])V,$\\
{\it and for other components, $\widetilde{R}$ equals to
$\overline{R}$.}\\

\noindent {\bf Proposition 6.4} {\it Let
$M=B\times_{b_1}F_1\times_{b_2}F_2\cdots
    \times_{b_m}F_m$ be a multiply twisted product and let $X,Y,Z\in \Gamma(TB)$
      and $V\in \Gamma(TF_i)$, $W\in \Gamma(TF_j)$, $U\in \Gamma(TF_k)$ and $P\in \Gamma(TF_l)$. Then}\\
\noindent
$(1)\widetilde{R}(V,W)X=\overline{R}(V,W)X+V(\pi(W))X-W(\pi(V))X-\pi([V,W])X.$\\
\noindent $(2)\widetilde{R}(U,V)W=U(\pi(V))W-V(\pi(U))W-\pi([U,V])W~ if ~i=k\neq j~ .$\\
\noindent $(3)\overline{R}(U,V)W=g(U,W){\rm grad}_B(V({\rm
ln}b_i))-g(V,W){\rm grad}_B(U({\rm
ln}b_i))+R^{F_i}(U,V)W-\frac{|{\rm
grad}_Bb_i|^2_B}{b_i^2}(g(V,W)U-g(U,W)V)
+g(W,\nabla_UP)V-g(W,\nabla_VP)U
+\pi(W)[\pi(V)U-\pi(U)V]+U(\pi(V))W-V(\pi(U))W-\pi([U,V])W
 ~if ~i=j=k=l,~$\\
{\it and for other components, $\widetilde{R}$ equals to
$\overline{R}$.}\\

By Proposition 6.3 and 6.4, we have\\

\noindent {\bf Proposition 6.5} {\it Let
$M=B\times_{b_1}F_1\times_{b_2}F_2\cdots
    \times_{b_m}F_m$ be a multiply twisted product and $P\in \Gamma(TB)$. Then $\widetilde{{\rm Ric}}=\overline{{\rm
    Ric}}.$}\\

    \noindent {\bf Remark.} By Proposition 6.5, then we can get the same
conclusions for $\widetilde{\nabla}$ and $\overline{\nabla}$.\\

\noindent {\bf Proposition 6.6} {\it Let
$M=B\times_{b_1}F_1\times_{b_2}F_2\cdots
    \times_{b_m}F_m$ be a multiply twisted product and $P\in \Gamma(TF_r)$. Then}\\
$$ \widetilde{{\rm Ric}}(V,W)=\overline{{\rm Ric}}
(V,W)+\sum_{j_r=1}^{l_r}\varepsilon_{j_r}g(V(\pi(E_{j_r}^r))W-E_{j_r}^r(\pi(V))W-\pi([V,E_{j_r}^r])W,E_{j_r}^r),\eqno(6.3)$$
{\it for $V,W\in\Gamma(TF_r)$. For other components,
$\widetilde{{\rm Ric}}=\overline{{\rm Ric}}.$}\\

 \noindent {\bf Acknowledgement.} This work
was supported by Fok Ying Tong Education
Foundation No. 121003.\\

\noindent{\large \bf References}\\

\noindent[AC1]N. Agashe, M. Chafle, A semi-symmetric non-metric
connection on a Riemannian manifold, Indian J. Pure Appl. math.
23(1992) 399-409.\\
\noindent[AC2]N. Agashe, M. Chafle, On submanifolds of a Riemannian
manifold with a semi-symmetric non-metric connection, Tensor(N.S.),
55(1994) 120-130.\\
 \noindent[ARS]L. Al\'{i}as, A.
Romero, M. S\'{a}nchez, Spacelike hypersurfaces of constant mean
curvature and Clabi-Bernstein type
problems, Tohoku Math. J. 49(1997) 337-345.\\
\noindent[BO]R. Bishop, B. O'Neill, Manifolds of negative curvature,
Trans. Am. Math. Soc. 145(1969) 1-49.\\
\noindent[DD]F. Dobarro, E. Dozo, Scalar curvature and warped
products of Riemannian manifolds, Trans. Am. Math. Soc. 303(1987)
161-168.\\
\noindent[DU]F. Dobarro, B. \"{U}nal, Curvature of multiply warped
products, J. Geom. Phys. 55(2005) 75-106.\\
\noindent[EJK]P. Ehrlich, Y. Jung, S. Kim, Constant scalar
curvatures on warped product manifolds, Tsukuba J. Math. 20(1996)
No.1 239-265.\\
\noindent[FGKU]M. Fern\'{a}ndez-L\'{o}pez, E. Garc\'{i}a-R\'{i}o, D.
Kupeli, B. \"{U}nal, A curvature condition for a twisted product to
be a warped product, Manu. math. 106(2001), 213-217.\\
\noindent[Ha]H. Hayden, Subspace of a space with torsion, Proc.
Lond. Math. Soc. 34(1932) 27-50.\\
 \noindent[SO1]S. Sular, C. \"{O}zgur, Warped products with a
semi-symmetric metric connection, Taiwanese J. Math. 15(2011) no.4
1701-1719.\\
\noindent[SO2]S. Sular, C. \"{O}zgur, Warped products with a
semi-symmetric non-metric connection, Arab J. Sci. Eng. 36(2011)
461-473.\\
\noindent[Wa]Y. Wang, Multiply twisted products, arXiv:1207.0199.\\
 \noindent[Ya]K. Yano, On semi-symmetric metric connection,
Rev. Roumaine Math. Pures Appl. 15(1970) 1579-1586.\\

\end{document}